\documentclass{article}
\usepackage[dvips]{color}
\usepackage{graphics,graphicx}
\usepackage{latexsym}
\usepackage{times}
\usepackage{cite}
\usepackage{natbib}  
\usepackage{amsbsy} 
\def\ts{\\[-12pt]}
\def\as{}
\def\ms{\medskip}

\def\ave{\mathop{\rm ave}}
\def\real{\mathop{{\rm I}\kern-.2em\hbox{\rm R}}\nolimits}
\def\footnoterule{\kern-3pt \hrule width 2truein \kern 2.6pt}
\def\dZ{\mathop{{\rm Z}\kern-.4em\hbox{\rm Z}}\nolimits}
\def\dC{\mathop{{\rm I}\kern-.5em\hbox{\rm C}}\nolimits}
\def\SUBSECTION#1\par{\vskip0pt plus.2\vsize\penalty-75
        \vskip0pt plus-.2\vsize\bigskip\bigskip
        \leftline{\bf #1}\nobreak\smallskip\noindent}
\def\leaderfill{\leaders\hbox to 1em{\hss.\hss.}\hfill}
\def\ba{{\bf a}}  \def\bz{{\bf z}}  \def\bw{{\bf w}}
\def\bx{{\bf x}}  \def\by{{\bf y}}    
\def\brho{\mathop{\rho\kern-.45em\hbox{$\rho$}}\nolimits}    %pm bold rho

  \def\b0{{\bf 0}}    
     \def\be{{\bf e}}  \def\bt{{\bf t}}
   
\def\bell{\mathop{\ell\kern-.47em\hbox{$\ell$}}\nolimits} %bold curly ell
\def\bep{\mathop{\epsilon\kern-.45em\hbox{$\epsilon$}}\nolimits}
\def\btheta{\mathop{\theta\kern-.45em\hbox{$\theta$}}\nolimits}
\def\bdelta{\mathop{\delta\kern-.45em\hbox{$\delta$}}\nolimits}
\def\blambda{\mathop{\lambda\kern-.55em\hbox{$\lambda$}}\nolimits}
\def\bnu{\mathop{\nu\kern-.50em\hbox{$\nu$}}\nolimits}
\def\bgamma{\mathop{\gamma\kern-.55em\hbox{$\gamma$}}\nolimits}
\def\bomega{\mathop{\omega\kern-.62em\hbox{$\omega$}}\nolimits}
\def\bxi{\mathop{\omega\kern-.50em\hbox{$\xi$}}\nolimits}
\def\bphi{\mathop{\phi\kern-.55em\hbox{$\phi$}}\nolimits}
\def\bmu{\mathop{\mu\kern-.55em\hbox{$\mu$}}\nolimits}
\def\btau{\mathop{\tau\kern-.50em\hbox{$\tau$}}\nolimits}
\def\balpha{\mathop{\alpha\kern-.54em\hbox{$\alpha$}}\nolimits}
\def\bbeta{\mathop{\beta\kern-.55em\hbox{$\beta$}}\nolimits}
\def\sqr#1#2{{\vcenter{\vbox{\hrule height.#2pt
               \hbox{\vrule width.#2pt height#1pt \kern#1pt
                \vrule width.#2pt}
              \hrule height.#2pt}}}}
 \def\square{\mathchoice\sqr34\sqr34\sqr{2.1}3\sqr{1.5}3}
\begin{document}
\baselineskip=18pt
\pagestyle{empty}
\centerline{\bf Breakdown Properties of the $M$-Estimators
of Multivariate Scatter}
\medskip

\centerline{David E.\ Tyler}
\centerline{Department of Statistics}
\centerline{Rutgers University}

\vskip 1cm
\begin{center}
\emph{Rutgers University Technical Report}, 1986 (unpublished).  \\
Abstract in: \emph{Inst.\ Math.\ Stat.\ Bull.}, 1986, Vol.\ 15, 116.
\end{center}

\vskip 1cm

\noindent
\centerline{\bf Summary} \ts

The $M$-estimates of multivariate scatter are known to have
breakdown points no greater than $1/(p+1)$, where $p$ is the
dimension of the data. In high dimension, the breakdown points
are usually considered to be disappointingly low. This paper studies
the breakdown problem in more detail. The exact breakdown points for the
$M$-estimates of scatter are obtained and it is shown that their low values
are primarily due to contamination restricted to some plane. If such
``coplanar'' contamination is not present, then there exists $M$-estimates
which have breakdown points close to $1/2$. The effect of ``coplanar''
contamination is further examined and is shown to be related to the
singularity of the scatter matrix. Finally, the implications of the
results of this paper on whether the low breakdown point is necessarily
a bad feature and on multivariate outlier detection are briefly discussed.

\vskip 2cm

\noindent
Research supported by NSF Grant DMS-8405325 

\newpage
\pagestyle{plain}
\noindent
{\bf 1. Introduction.}

The affine invariant $M$-estimates of multivariate location and scatter
were first proposed by Maronna (1976) as robust alternatives to the
sample mean vector and covariance matrix. One feature of these estimators,
though, which was noted by Maronna (1976) and has been a concern to others,
{\it e.g.}, Huber (1981), Stahel (1981), Donoho (1982), and Devlin {\it et al.}
(1981), is their relatively low breakdown point, particularly in higher
dimensions. Maronna (1976) obtains an upper bound for the breakdown point
of an $M$-estimator and shows that none have a breakdown point greater
than $1/(p+1)$, where $p$ represents the dimension of the data.
Stahel (1981) obtains a general bound of $1/p$ for a slightly more
general class of $M$-estimators. Although much work has appeared on
properties and applications of the $M$-estimators of multivariate
location and scatter, there has been no further theoretical results on
their breakdown properties.

The aim of this paper is to study the breakdown problem in more detail and
to address the question: Is the low breakdown point necessarily a bad
feature?  The notion of breakdown is viewed here more as a descriptive
rather than an optimal property. Attention is restricted to the $M$-estimates
of scatter in this paper since the low breakdown point of the multivariate
$M$-estimates is due to the breakdown of the scatter component, as 
demonstrated by both Maronna (1976) and Stahel (1981).

Loosely summarizing, it is shown in Section 3 that the upper bounds given by
Maronna (1976) for the breakdown point of the $M$-estimates of scatter are
in fact the exact breakdown points. In Section 4, the cause of the low
breakdown point is investigated and is shown to be primarily due to
contamination restricted to some plane, a type of contamination unique to
the multivariate setting. In fact, if ``coplanar'' contamination is not
present, then there exist $M$-estimates with breakdown points close
to 1/2 (Theorem 4.1). Furthermore, some $M$-estimates of scatter are
shown to breakdown under a small percent of ``coplanar'' contamination,
even though no ``outliers'' or ``inliers'' are present (Theorem 4.2).
Section 5 examines the effect of ``coplanar'' contamination, which as
one might expect, is related to the singularity of the scatter matrix.

After formally presenting the aforementioned results, some brief concluding
remarks concerning their implications are made in Section 6. To begin, some
background on the $M$-estimates of scatter and on finite sample breakdown
is given.
\bigskip

\noindent
{\bf 2. Background.}

{\it 2.1 $M$-estimators of scatter.}
For $p$-dimensional data $x_1,x_2,\ldots,x_n$,  Maronna (1976) defines the
affine invariant $M$-estimator of scatter about some fixed center $\bt$
to be the positive definite symmetric (p.d.s.) matrix $V_n$ satisfying
the equation
$$V_n = \ave \big\{u(s_i)(\bx_i - \bt)(\bx_i - \bt)'\big\}\leqno(2.1)$$
where $s_i = (\bx_i - \bt)' V_n^{-1}(\bx_i  - \bt)$ and $u$ is some
scalar valued function. The $M$-estimator $V_n$ can be viewed as an
adaptively weighted covariance matrix whose weights depend on an
adaptive Mahalanobis distance from the center. For future reference,
multiplying (2.1) by $V_n^{-1}$ and taking the trace gives
$$p = \ave \big\{\psi(s_i)\big\},\leqno(2.2)$$
where $\psi(s) = su(s)$. Also, let  $K = \sup_{s > 0} \psi(s)$.

Some conditions on the function $u$ and on the empirical distribution are
needed to insure the existence and uniqueness of $V_n$. The existence
lemma given below is from Tyler (1985). \ms

C{\footnotesize{ONDITION}} 2.1.

(i)  $u(s)$ is non-negative, non-increasing and continuous for $s > 0$.

(ii) $u(s)$ and $s$ are bounded.

(iii) $\psi(s)$ is non-decreasing for $x > 0$ and strictly increasing for
$\psi(s) < K$.

(iv)  $K > p$. \ms

\noindent
Let $n_0$ represent the size of the largest subset of
$X = \{ \bx_1,\bx_2,\ldots,\bx_n\}$ which is in general position about
the center $\bt$, where a set of vectors from $\real^p$ is said to be in
general position about a fixed vector $\bt$ if the plane generated by any
subset of size $p$ together with $\bt$ is $\real^p$. Let $P_n$ be the
empirical distribution function of \mbox{$\big\{ (\bx_i - \bt); 1\le i\le n\big\}$.}  \ms

C{\footnotesize{ONDITION}} 2.2.  For any subspace $S$ with $0\le {\rm rank} (S)
\le m-1$,

(i) $P_n(S) < 1 - p/K + \min[1, n_o\, {\rm rank} (S)/n]/K$ and 
$n_0 > p(p-1)$.

(ii)  $P_n(S) \le 1 - \{ p = {\rm rank} (S)\}/K$.  \ms

\newpage
L{\footnotesize{EMMA}} 2.1. {\it
Suppose $u$ satisfies Condition 2.1.

(i)  If Condition 2.2.i holds, then there exists a unique p.d.s.\ solution
$V_n$ to (2.1).

(ii)  If a p.d.s.\ solution $V_n$ exists to (2.1), then Condition 2.2.ii
must hold.

(iii)  If a p.d.s.\ solution $V_n$ exists to (2.1) and $n_0 > p$,
then it is unique.} \ms

Lemma 2.1 essentially states nonexistence of $V_n$ is due to too many
data points being coplanar with the center $\bt$.  

Maronna (1976) and Huber (1981) also give sufficient conditions for
existence of $V_n$. Huber's condition on $u$ is more general than
Condition 2.1. Both Huber's and Maronna's  condition of $P_n$ are
more restrictive than Condition 2.1.i. \ms

{\it 2.2 Finite sample breakdown.}
A number of different definitions of the breakdown point of an estimator have
been proposed since Hampel (1971) formally introduced the concept. Recently,
Donoho (1982) and Donoho and Huber (1983) define the notation of finite
sample breakdown in the following manner. Let $m$ arbitrary data points
$Y = \{\by_1,\by_2,\ldots,\by_m\}$ augment the original data
$X = \{\bx_1,\bx_2,\ldots,\bx_n\}$ producing an $\epsilon$-contaminated
sample $Z = X \cup Y$ consisting of a fraction of $\epsilon = 
m/(n+m)$ bad values. For a given $\epsilon$, a statistic is said to
breakdown under $\epsilon$-contamination if the difference between the
statistic defined on the original sample $X$ and the statistic defined
on the contaminated sample $Z$ can be made arbitrarily large in some sense
for varying choices of $Y$. The finite sample breakdown point of the
statistic at the sample $X$ is $\epsilon^*(X)$, the infimum of all
$\epsilon$ producing breakdown.

Let $V_n(X)$ and $V_{n+m}(Z)$ represent p.d.s.\ solutions to (2.1) for
the original data $X$ and the contaminated data $Z$ respectively whenever
they exist. For $\epsilon = m/(n+m)$ and $V_n(X)$ existing, define the
maximum ``bias'' at $X$ caused by $\epsilon$-contamination to be \\[-8pt]
\[b(\epsilon; X) = \cases{\sup [{\rm trace} \{V_{n+m}(Z)V_n^{-1}(X)
                                 + V_n(X)V_{n+m}^{-1}(Z)\}],&
                                      $Z\in S_m(X)$,\cr
                           \infty,&$Z{\not\in} S_m(X)$,\cr}\]			
\\[-8pt]																		
where $S_m(X) = \{Z = X\cup Y| V_{n+m}(Z)\,\, {\rm exists}\}$ and the
supremum is taken over all choices of $Y$ and all possible solutions for
$V_n(X)$ and $V_{n+m}(Z)$.  Breakdown occurs under $\epsilon$-contamination
whenever $b(\epsilon; X) = \infty$. This implies either the statistic
$V_{n+m}(Z)$ does not exist, trace$\{V_{n+m}(Z)\}$ can be made 
arbitrarily large or $V_{n+m}(Z)$ can be made arbitrarily close to the
zero matrix or some other singular matrix. This notion of breakdown for a
p.d.s.\ statistic is in agreement with the notion used by Maronna (1976),
Stahel (1981), and Donoho (1982). The finite sample breakdown point of
$V_n(X)$ at $X$ is defined to be
$$\epsilon^*(X) 
     = \min_m \big\{\epsilon = m/(n+m)| b(\epsilon, X) = \infty\big\}.$$

To simplify notation, the results of this paper are stated in terms of
$\delta^*(X)$ where $\delta^*(X)$ is defined to be a fraction such that
$b(\epsilon, X) = \infty$ if $\epsilon > \delta^*(X)$ and
$b(\epsilon, X) < \infty$ if $\epsilon < \delta^*(X)$.
Since the possible values of $\epsilon$ are discrete, $\delta^*(X)$ is not
uniquely defined. The relationship between $\delta^*(X)$ and $\epsilon^*(X)$
is easily shown to be
$$\delta^*(X) \le \epsilon^*(X)
               < \big[\delta^*(X) + \{\delta^*(X)\}/n\big]/
                 \big[1 + \{1 - \delta^*(X)\}/n].\leqno(2.5)$$

\bigskip
\noindent
{\bf 3. The Breakdown Point of $V_n(X)$.}

Hereafter, assume that the ``good'' data $X$ is in general position about
$\bt$, which occurs almost surely when sampling from a continuous distribution
in $\real^p$. This assumption concerning $X$ is also used by Donoho (1982) in
studying the finite sample breakdown properties of projection pursuit based
estimators of location and scatter. It is also assumed hereafter that
$n > p(p-1)$. By Lemmas 2.1.i and 2.1.ii, these assumptions assure the
existence and uniqueness of $V_n(X)$ and the uniqueness of $V_{n+m}(Z)$
if it exists.

The general breakdown point of $V_n(X)$ is given in Theorem 3.1 below.
Before presenting the theorem some lemmas concerning the existence and
behavior of $V_{n+m}(Z)$ are given. The proofs of the lemmas are given
in the appendix. For brevity, let $\epsilon_m = m/(n+m)$, and in all
proofs assume without loss of generality that $\bt = {\bf 0}$. \ts

L{\footnotesize{EMMA}} 3.1.  {\it $Z\in S_m(X)$ if either

(i)  $\epsilon_m < 1 - p/K$ and $n \ge K$,

(ii)  $\epsilon_m < 1 - np/\{nK - (K-n)(p-1)\}$ and $n < K$, or

(iii)  $\epsilon_m < 1 - \max [np/\{nK + (n-K)\}, n(p-1)/\{(n-p+1)K\}]$,
\mbox{$n > K$ and $\bt {\not\in} Z$.}
} 

\newpage
L{\footnotesize{EMMA}} 3.2. {\it

(i)  If $\epsilon_m < 1 - p/K$, then $\{ {\rm trace}\, V_{n+m}^{-1}(Z)|
Z \in S_m(X)\}$ is bounded above.

(ii)  If $\epsilon_m < 1/K$, then $\{ {\rm trace}\, V_{n+m}(Z)|
Z \in S_m(X)\}$ is bounded above.

(iii)  If $\epsilon_m < p/K$, then $\{ {\rm trace}\, V_{n+m}^{-1}(Z)|
Z \in S_m(X)\}$ is bounded away from zero.} \ts

T{\footnotesize{HEOREM}} 3.1. {\it 
$\delta^*(X) = \min \{1/K, 1 - p/K\}$ for $n+1 > K$, and
$\delta^*(X) = 0$ for $n+1 \le K$.} \ts

{\bf Proof:}
If $Y = \{ {\bf 0},{\bf 0},\ldots,{\bf 0} \}$ and $\epsilon_m > 1 - p/K$,
then $P_{n+m}({\bf 0}) = \epsilon_m > 1 - p/K$.
This implies by Lemma 2.1.ii that $V_{n+m}(Z)$ does not exist.
If $Y = \{\by, \by,\ldots,\by\}$ and $V_{n+m}(Z)$ exists, then
$\by' A\by = \ave \{u(\bz' A\bz)(\by' A\bz)^2\}$ where
$A = \{V_{n+m}(Z)\}^{-1}$ and the average is over $z \in Z$.
This implies $ 1 - \ave \{u(\bz' A\bz)(\by' A z)^2/\bz' A\bz\}$ or
$$1 = \epsilon_m \psi(\by' A\by) + %\break
 (n+m)^{-1} 
     \sum_{1\le i\le n} u(\bx'_i A\bx_i)
     \times (\bx'_i A\by)^2/\by' A\by.$$
Express $\by = r {\btheta}$ where $\btheta'\btheta = 1$ and let $r\to\infty$.
If $V_{n+m}(Z)$ does not breakdown as $r\to\infty$, the $\psi(\by' A\by)
\to K$ and since $X$ spans $\real^p$, for some $\bx \in X$, 
$u(\bx' A\bx)(\bx'A\by)^2/\by' A\by 
= u(\bx' A\bx)(\bx' A\btheta)^2/\btheta' A\btheta$ does not go to zero.
This implies $\epsilon_m < 1/K$ and when using $m = 1$,
$n+1 > K$. Thus $\delta^*(X) \le \min (1 - p/K, 1/K)$ and if $n+1 \le K$,
$\delta^*(X) = 0$.

If $\epsilon_m < \min (1/K, 1 - p/K)$, then by Lemma 3.1, $z\in S_m(X)$.
Application of Lemma 3.2 gives $\delta^*(X) \ge \min (1/K, 1 - p/K)$.
$\qquad\square$ \ms

Maronna (1976) obtains $\min \{1/K, 1 - p/K\}$ as an upper bound of the
$M$-estimator of scatter at any continuous elliptically contoured
distribution in $\real^p$, and conjectures that the bound is the exact
breakdown point. He uses the definition of breakdown at a model
distribution rather than finite sample breakdown. The arguments
given in the proof of Lemma 3.2 and Theorem 3.1 can be modified
to show that for the case of known center the breakdown point of the
$M$-estimate of scatter is equal to $\min \{1/K, 1 - p/K\}$ at any
continuous model in $\real^p$. Maronna further states that this
upper bound is obtained by letting a point mass contamination go
to infinity. This is true for the $1/K$ term but not the $1 - p/K$
term, which is obtained by point mass contamination at the center.
As noted by Maronna, the breakdown point is low for higher dimensions
since $K > p$ and so $\min \{1/K, 1 - p/K\} \le 1/(p+1)$.
\bigskip

\noindent
{\bf 4. The Sources of Breakdown.}

The objective of this section is to investigate what causes the
$M$-estimate of multivariate scatter to breakdown. For univariate
scale problems, breakdown is usually due to the existence of too
many outliers or to the existence of too many inliers, that is, data
points near the center. In the multivariate setting, though, breakdown
may also occur because of too many data points lying in some lower
dimensional plane containing the center of $\bt$, which will be 
referred to as coplanar contamination. By examining the proof of
Theorem 3.1, one can note that the low overall breakdown point of the
$M$-estimates of scatter, that is the $1/K$ term, is obtained by
outliers which are coplanar with the center. If coplanar contamination
is not present, then it is shown in Theorem 4.1 below that some
$M$-estimators of scatter can have breakdown points close to 1/2.
Before formally presenting this result, some additional notation
and definitions are needed.

Let $C_m(X)$ be a subset of the product set
$\prod_{j=1}^m \real^p$, possibly dependent on $X$. Define the finite sample
breakdown point of $V_n$ at $X$ due to a sequence $C_1(X),C_2(X),\ldots$
to be $\epsilon^*(X)$ where $\epsilon^*(X)$ is defined by (2.4) but with
the restriction $Y \in C_m(X)$ in the definition of $b(\epsilon, X)$.

An element $\bz \ne \bt$ from $\real^p$ can be expressed as $\bz = \bt + r\btheta$ where $\theta = (\bz - \bt)/r$ and
$r = \{(\bz - \bt)'(\bz - \bt)\}^{1/2}$. Using this representation, define for
$Z = X \cup Y = \{\bz_1,\bz_2,\ldots,\bz_{n+m}\}$
$$\rho_m(Z) 
= \min \lambda_p\bigg\{\sum_{j=1}^p \btheta_{i(j)} \btheta'_{i(j)}\bigg\}
\leqno(4.1)$$
where the minimum is taken over all subsets of size $p$ from $Z$ for
which $\bz\ne t$, and $\lambda_p(\cdot)$ represents the smallest
eigenvalue of the $p\times p$ non-negative definite argument.
The quantity $\rho_m(Z) \ne 0$ if and only if
$\{\bz| \bz\in Z, \bz\ne t\}$ is in general position about $\bt$.
Also, define for $Y = \{\by_1,\by_2,\ldots,\by_m\}$
$$r_m(Y) = \min\{ (\by_i - \bt)'(\by_i - \bt),\,\, 1\le i\le m\},
\leqno(4.2)$$
and let $C_{1,\rho,m}(X) = \{Y| \rho_m(X\cup Y) > \rho\}$,
$C_{2,r,m}(X) = \{Y| r_m(Y) > r\}$ and
$C_{3,r,m} = \{Y|r_m(Y) < B\}$.

Some results concerning the behavior of $V_{n+m}(Z)$ when $Y$ is
restricted to certain classes are given in the following lemma.
The proof of the lemma is given in the appendix. \ms

L{\footnotesize{EMMA}} 4.1. {\it Let $\rho > 0,\,\, r > 0$ and $B < \infty$.

(i)  If $n+1 > K$ and $\epsilon_m < p/K$, 
then $\{ {\rm trace} \,\, V_{n+m}(Z)| Z\in S_m(X)\,\, and\,\,
Y\in C_{1,p,m}(X)\}$ is bounded above.

(ii) For any $m$, $\{ {\rm trace}\,\, V_{n+m}(Z)| Z\in S_m(X)\,\, and\,\,
Y\in C_{3,B,m}(X)\}$ is bounded above.

(iii) For any $m$, $\{ {\rm trace}\,\, V_{n+m}(Z)| Z\in S_m(X)\,\, and\,\,
Y\in C_{2,r,m}(X)\}$ is bounded away from zero.

(iv) If $\epsilon_m < 1 - n(p-1)/\{(n-1)K\}$, then
$\{ {\rm trace}\,\, V_{n+m}^{-1}(Z)| Z\in S_m(X)\,\, 
and\,\,Y\in C_{2,r,m}(X)\}$ is bounded above.

(v) If $n+1 > K$, then for any $m$,
$\{ {\rm trace}\,\, V_{n+m}^{-1}(Z)| Z\in S_m(X)\,\, and \,\,
Y\in C_{1,p,m}(X) \cap C_{2,r,m}(X)\}$ is bounded above.} \ms

The following results concerning the breakdown of $V_n(X)$ whenever
coplanar contamination is not present are similar to the breakdown
results for univariate scale. Estimators which protect against
outliers, {\it i.e.} $K$ near $p$, tend to breakdown in the
presence of inliers. For the compromising choice $K = 2p$, the 
breakdown point given in Theorem 4.1.iii is approximately 1/2. \ms

T{\footnotesize{HEOREM}} 4.1. {\it Let $\rho > 0,\,\, r > 0$ and $B < \infty$.

(i)  For the sequence $C_{1,\rho,m}(x) \cap C_{2,r,m}(X)$,
$\delta^*(X) = p/K$ if $n+1 > K$ and $\delta^*(X) = 0$ if
$n+1 \le K$.

(ii)  For the sequence $C_{2,B,n}(X)$,
$\delta^*(X) = 1 - p/K$ for $n \ge K$ and
$1 - np/\{nK - (K-n)(p-1)\} \le \delta^*(X) \le 1 - p/K$ for $n < K$.

(iii)  For the sequence $C_{1,\rho,m}(X)$,
$\delta^*(X) = \min (p/K, 1 - p/K)$ if $(n+1) > K$ and 
$\delta^*(X) = 0$ if $n+1 \le K$.
} \ms

{\bf Proof:} 
(i)  If $y\in C_{1,\rho,m}(X) \cap C_{2,r,m}(X)$, 
then it can be verified from Lemma 2.1.i that $Z\in S_m(X)$ for $n > p(p-1)$
after noting $P_{n+m}(S) \le {\rm rank} (S)/(n+m)$ and
$(n+m)_0 = n+m$.  For $n+1 > K$, it follows from Lemmas 4.1.i and 4.1.v that
$\delta^*(X) \ge p/K$.  For (2.2), $p \ge 
(n+m)^{-1} \sum_{1\le i\le m} \psi(\by'_i V_{n+m}^{-1}(Z)\by_i) \to\epsilon_m$
if $V_{n+m}(Z)$ does not breakdown as $\by'_i \by_i \to\infty$,
$1\le i\le m$, and so $\delta^*(X) \le p/K$. For $n+1 \le K$, the proof that
$\delta^*(X) = 0$ is analogous to the proof in Theorem 3.1.

(ii)  As in the proof of Theorem 3.1, the upper bound for $\delta^*(X)$ is
obtained by choosing $Y = \{ {\bf 0},{\bf 0},\ldots,{\bf 0}\}$. The lower
bound follows from Lemmas 3.1, 3.2.i and 4.1.ii.

(iii) The lower bound follows from Lemmas 3.1.i, 3.2.i and 4.1.i. The upper
bound follows from parts (i) and (ii) of this theorem.$\qquad\square$ \ms

Theorem 4.1.i generalizes a statement made by Maronna (1976) in which he
quotes $p/K$ as the breakdown point due to contaminating a spherically
contoured model distribution by a long-tailed spherically contoured
distribution.

An interesting aspect to the multivariate breakdown problem is that
breakdown can occur because of coplanar contamination, even though
the contamination contains no outliers or inliers. In fact, as seen
in the next theorem, the breakdown point due  to such contamination
can be quite low. \ms

T{\footnotesize{HEOREM}} 4.2. {\it 
Let 
$a_1 = 1 - np/\{nK - (K-n)(p-1)\}$,
$a_2 = 1 - np/(nK + n-K)$,
$a_3 = 1 - n(p-1)/\{(n - p+1)K\}$, and
$a_4 = 1 - n(p-1)/\{(n-1)K\}$.
For the sequence $C_{2,r,m}(X) \cap C_{3,B,m}(X)$, with $r > 0$ and
$B < \infty$, $\min (a_2,a_3) \le \delta^*(X) \le a_4$ when $n > K$, and
$a_1 \le \delta^*(X) \le a_4$ when $n\le K$.} \ms

{\bf Proof:}
The upper bound is obtained by letting $Y = \{\bx_1,\bx_1,\ldots,\bx_1\}$
and applying Lemma 2.1.ii.  The lower bounds follow from Lemmas 3.1.ii,
3.1.iii, 4.1.ii and 4.1.iv.$\qquad\square$ \ms

As $n$ goes to infinity, the bounds on $\delta^*(X)$ in Theorem 4.2
simplify to
$$1 - p/(K+1) \le \delta^*(X) \le 1 - (p-1)/K.\leqno(4.3)$$
For $K$ near $p$, the breakdown point is in the neighborhood of
$1/(p+1)$ to $1/p$.

It is interesting to note that for $K$ near $p$, the influence  of
non-coplanar outliers is essentially nonexistent; see Theorem 4.1.i.
Furthermore, for such $K$ the breakdown points given in Theorem 4.2
do not differ greatly from $1/K$. Therefore, for such
$M$-estimators of scatter the low breakdown point caused by coplanar
outliers can be attributed primarily to the coplanar aspect of the
contamination rather than the outlier aspect.  A brief heuristic
explanation is helpful in understanding this phenomena. The defining
equation (2.1) can be rewritten as
$$V_n = n^{-1} \sum\psi (s_i)(\bx_i - \bt)(\bx_i - t)'/s_i\leqno(4.4)$$
where the summation is over $\bx_i \ne \bt$. The function $\psi$ can be
viewed as measuring the influence of the distance of an observation from
the center. The term $(\bx_i - \bt)(\bx_i - \bt)'/s_i$ is dependent only
on the direction of the observation from the center and not on the
distance. Since $\psi$ is non-decreasing, if $K$ is near $p$, then (2.2)
implies that $\psi$ is roughly a constant function. Thus, outliers have
little more influence than other data points and breakdown is primarily
dependent on the interrelationships of the directions of the data points
from the center.
\bigskip

\noindent
{\bf 5. The Effect of Coplanar Contamination.}

The notion of contamination which is coplanar with the center
distinguishes the multivariate breakdown problem from the univariate one.
Intuitively, one might expect such contamination would be related in some
way to the singularity of the estimate of scatter. In this section, this
intuition is briefly but formally investigated.

The difference in the breakdown points in Theorem 3.1 and Theorem 3.2.iii
can be attributed to the existence of outliers which are coplanar with the
center, and by Theorem 4.2 cannot be attributed to coplanar contamination
alone. For $n+1 > K$, and
$$\min (1/K, 1 - p/K) < \epsilon_m < \min (p/K, 1 - p/K),\leqno(5.1)$$
Lemma 3.1.ii implies $Z\in S_m(X)$, and furthermore Lemmas 3.2.i and 3.2.iii
imply that there exists a nonzero non-negative definite symmetric matrix
$A_1$ and a positive definite symmetric matrix $A_2$ such that for all
$Z \in S_m(X)$
$$A_1 < V_{n+m}^{-1}(Z) < A_2,\leqno(5.2)$$
where the ordering refers to the partial ordering of symmetric matrices.
Thus, for $n+1 > K$ and (4.2) holding, ``coplanar outliers'' tend to make
$V_{n+m}^{-1}(Z)$ singular rather than ``blowing up'' or becoming strictly
zero. A natural question which arises is whether the limiting null space
of $V_{n+m}^{-1}(Z)$ and the contaminating plane coincide. For the
following case, which produces the $1/K$ term in Theorem 3.1, they do.
The proof is given in the appendix. \ms

T{\footnotesize{HEOREM}} 5.1. {\it  For fixed $m$, let $Z_r = X \cap Y_r$
where $Y_r = \{\by_r,\by_r,\ldots,\by_r\}$, $\by_r = r\btheta + \bt$
and $\btheta'\btheta = 1$. If $n+1 > K$ and (5.2) holds, then
$\btheta' V_{n+m}^{-1}(Z_r)\btheta \to 0$ as $|r| \to\infty$, and if
$\ba$ is not proportional to $\btheta$, then 
$\inf_r \ba'V_{n+m}^{-1}(Z_r)\ba >0$.} \ms

As shown in Theorem 4.2, contamination within some hyperplane containing
the center $\bt$ can cause breakdown, even though no outliers or inliers
are present. Breakdown by such contamination is due to either nonexistence
or to the $M$-estimate of scatter tending toward singularity. However,
the estimate does not tend to zero nor does it become arbitrarily large.
To state this formally, let $r > 0$ and $B < \infty$. Lemmas 4.1.ii and 4.1.iii
imply there exists a nonzero non-negative definite symmetric matrix $W_1$ and
a positive definite symmetric matrix $W_2$ such that for all $Y\in
C_{2,r,m}(X) \cap C_{3,B,m}(X)$
$$W_1 < V_{n+m}(Z) < W_2\leqno(5.3)$$
provided $Z = X\cup Y\in S_m(X)$. Furthermore, if one considers a sequence
$Y_k\in C_{2,r,m}(X) \cap C_{3,B,m}(X)$ such that
$Z_k = X\cup Y_k\in S_m(X)$ and $Z_k\to Z {\not\in} S_m(X)$, then (5.3)
and the continuity of $u$ imply that the largest root of
$V_{n+m}(Z_k)$ is bounded away from zero and infinity for all $k$, and the
smallest root tends to zero. Again, a natural question which arises is
whether the limiting range of $V_{n+m}(Z_k)$ and the contaminating plane
coincide. For the following case, which produces the upper bound in
Theorem 3.3, they do. The condition $\epsilon_m < 1 - n(p-2)/\{K(n-2)\}$
is probably not needed in the following theorem, but the author is not
able to derive the result without this condition. The proof is given
in the appendix. \ms

T{\footnotesize{HEOREM}} 5.2. {\it For fixed $m$, let $Z_k = X \cup Y_k$ where
$Y_k = \{\bx_1 + \bw_{1,k}, \bx_1 + \bw_{2,k},\ldots,\bx_1 + \bw_{m,k}\}$
with $\bw_{i,k} \to {\bf 0}$, $1\le i\le m$, and $Z_k\in S_m(X)$ for all $k$.
If $$1 - n(p-1)/\{K(n-1)\} < \epsilon_m < 1 - n(p-2)/\{K(n-2)\},$$ then
$\inf_k \bx'_1 V_{n+m}(Z_k)\bx_1 > 0$, and for any $\ba$ such that
$\ba'\bx_1 = 0$, $\ba'V_{n+m}(Z_k)\ba \to 0$.}
\bigskip
\newpage

\noindent
{\bf 6. Concluding Remarks.}

Is the low breakdown point necessarily a weak feature of an
$M$-estimate of multivariate scatter? One can respond yes if it is
believed that contamination lying in or near some lower dimensional 
plane is feasible and no attempt is made to detect such contamination.
Otherwise, $M$-estimators exist which have good breakdown properties.

An alternative or complimentary approach is to try to detect bad data points,
particularly outliers. The results of Section 6 suggest that if a group of
outliers lie in or near some lower dimensional plane, then the near
singularity of an $M$-estimate of scatter can be used to help detect such
systematic outliers, with the directions associated with the largest roots
indicating where to search for the outliers. More research along these lines
may be fruitful. If outliers exist which are not coplanar, then their
detection may be more difficult. $M$-estimates of scatter exist, though,
which are quite stable under such contamination.

Finally, if coplanar contamination is present, with or without outliers,
it may be desirable to note this rather than simply attempt to
summarize the data via a location and scatter statistics. Again the
results of Section 6 suggest the near singularity of an $M$-estimate
of scatter may indicate the existence of such systematic contamination,
with the directions associated with the larger roots coinciding with the
contaminating plane.
\bigskip

\noindent
{\bf 7. Appendix:  Some Proofs.}

The proofs of Lemmas 3.1, 3.2, 3.3 and 4.1 and Theorems 5.1 and 5.2 are given
in this appendix. In all proofs, without loss of generality, $\bt$ is set
equal to {\bf 0}.  Recall $X$ is assumed to be in general position about
$\bt = {\bf 0}$, and $n > p(p-1)$. \as

{\bf Proof of Lemma 3.1:}
This lemma follows from Lemma 2.1.i after noting that
$P_{n+m}(S) \le \{m + {\rm rank} (S)\}/(n+m)_0 \ge n$. Further,
if $\bt {\not\in} Z$, then $P_{n+m}({\bf 0}) = 0$. These above
statements are true since $X$ is in general position about the
center.$\qquad\square$ \as

{\bf Proof of Lemma 3.2:}
Consider any sequence $V_{n+m}(Z_k)$ where $Z_k = X\cup Y_k \in S_m(X)$.
Let $\alpha_k = {\rm trace}\, \{V_{n+m}^{-1}(Z_k)\}$ and
$\Gamma_k = V_{n+m}^{-1}(Z_k)/\alpha_k$. Since ${\rm trace}\, \Gamma_k = 1$,
there exists a convergent sub-sequence, say for $j\in J$,
$\Gamma_j \to\Gamma$ a positive semi-definite symmetric (p.s.d.s.) matrix
with ${\rm trace}\, (\Gamma) = 1$. Let $Y_k = \{\by_{i,k}; 1\le i\le m\}$,
$r_{i,k} = (\by'_{i,k} \by_{i,k})^{1/2}$ and for $r_{i,k}\ne 0$, let
$\btheta_{i,k} = \by_{i,k}/r_{i,k}$. If $r_{i,k} = 0$, define
$\btheta_{i,k} = \be$, a vector such that $\be'\be = 1$ and
$\Gamma\be \ne 0$. Such a vector exists since ${\rm trace}\, \Gamma = 1$
and so $\Gamma\ne 0$. Since $\btheta'_{i,k}\btheta_{i,k} = 1$, the
sub-sequence $J$ can be chosen so that for each $1\le i\le m$,
$\btheta_{i,j} \to\btheta_i$ with $\btheta'_i\btheta_i = 1$.
Pre- and post-multiply (2.1) by $\alpha_k^{1/2}\Gamma_k^{1/2}$,
where $\Gamma_k^{1/2}$ is the unique p.s.d.s. square root of
$\Gamma_k$, and then multiply by the orthogonal projection into the null
space of $\Gamma$, say $P_\Gamma$. This gives
$$P_\Gamma = \ave \big\{u(\alpha_k\bz'\Gamma_k\bz)
                        \alpha_k P_\Gamma \Gamma_k^{1/2}\bz\bz'\Gamma_k^{1/2} 
                          P_\Gamma \big\}\leqno(7.1)$$
where the average is over $z\in Z_k$. Taking the trace gives 
${\rm rank}\,(P_\Gamma) = \ave \big\{\psi_k(\bz)\big\}$, where
$\psi_k(\bz) = u(\alpha_k \bz'\Gamma_k \bz)\alpha_k \bz'\Gamma_k^{1/2}
P_\Gamma \Gamma_k^{1/2} \bz.$ If $\Gamma\bx_i \ne {\bf 0}$, then
$$\psi_k(\bx_i) = \psi(\alpha_j\bx'_i\Gamma_j\bx_i)\times
\bx'_i\Gamma_j^{1/2} P_\Gamma \Gamma_j^{1/2}\bx_i/\bx'_i\Gamma_j\bx_i,$$ which
goes to zero since $\psi$ is bounded and $P_\Gamma \Gamma_j^{1/2} = 0$.
Likewise, if $\Gamma\btheta_i \ne {\bf 0}$, then
$\psi_j(\by_{i,j})
= \psi(\alpha_j\by'_{i,j} \Gamma_j\by_{i,j})
\btheta'_{i,j} \Gamma_j^{1/2} P_\Gamma \Gamma_j^{1/2}\- \btheta_{i,j}/
\btheta'_{i,j} \Gamma_j \btheta_{i,j}
\to 0$, and so 
$$\lim_{j\in J} (n+m)^{-1} \sum_{\bz\in Z_{0,j}} \psi_j(\bz) = R_\Gamma,
\leqno(7.2)$$
where $Z_{0,j} = \{\bx_i| \Gamma\bx_i = {\bf 0}\}
                 \cup \{\by_{i,j}| \Gamma\btheta_i = {\bf 0}\}$
and $R_\Gamma = {\rm rank}\, (P_\Gamma)$.

From (2.2), 
$p \ge (n+m)^{-1}\big\{\sum_{z\in Z_{0,j}} \psi(\alpha_j \bz'\Gamma_j\bz)
                        + \sum_{\Gamma\bx_i \ne 0}
                                  \psi(\alpha_j\bx'_i\Gamma_j\bx_i)\big\}$.
If $\Gamma\bx_i \ne 0$ and $\alpha_j\to\infty$ for $j\in J$, then
$\psi(\alpha_j\bx'_i\Gamma_j\bx_i)\to K$. Since $P_\Gamma$ is idempotent,
$\psi_j(\bz) \le u(\alpha_j\bz'\Gamma_j\bz)\alpha_j\bz'\Gamma_k\bz 
= \psi(\alpha_j\bz'\Gamma_j\bz)$. These results, together with (7.2) and
the assumption that $X$ is in general position about the origin, imply
that if $\alpha_j \to\infty$ for $j\in J$, then
$p \ge R_\Gamma + (n-R_\Gamma)K/(n+m)$. This last inequality is equivalent to
$$\epsilon_m \ge 1 - nn(p-R_\Gamma)/\{K(n-R_\Gamma)\}.\leqno(7.3)$$

(i) The right-hand side of (7.3) is an increasing function or $R_\Gamma$
for $0\le R_\Gamma \le p$. Thus, if $\epsilon_m < 1 - p/K$, then $\alpha_k$
must be bounded above.

(ii) If $\alpha_k$ is not bounded above, then $J$ can be chosen so that
$\alpha_j \to\infty$ for $j\in J$. The right-hand side of (7.3) is greater
than $1/K$  unless $R_\Gamma = 0$ and thus $\Gamma$ is nonsingular. This
implies ${\rm trace}\, \{V_{n+m}(Z_j)\} = {\rm trace}\, \Gamma_j^{-1}/\alpha_j
\to 0$.

If $\alpha_j$ is bounded above, then $J$ can be chosen so that
$\alpha_j\to\alpha$ for $j\in J$. This implies for $\Gamma\bx_i = {\bf 0}$,
$\psi_j(\bx_i) \le \psi(\alpha_j\bx'_i\Gamma_j\bx_i)\to \psi(0) = 0$ and so
from (7.2), ${\rm rank}\, (P_\Gamma) \le \epsilon_m K$. This contradicts the
condition on $\epsilon_m$ unless ${\rm rank}\, (P_\Gamma) = 0$ and thus
$\Gamma$ is nonsingular. This implies ${\rm trace}\, \{V_{n+m}(Z_j)\} \to 
{\rm trace}\, \Gamma^{-1}/\alpha$ unless $\alpha = 0$. If $\alpha = 0$, then
$\psi(\alpha_j\bx'_j\Gamma\bx_j) \to 0$ and so by (2.2), $p\le\epsilon_m K$,
which contradicts the condition on $\epsilon_m$.

(iii) If $\alpha_k$ is not bounded away from zero, the $J$ can be chosen so
that $\alpha_j \to 0$ for $j\in J$. Using (7.2), this implies $R_\Gamma =0$
and so $\Gamma$ is nonsingular. By (2.2), this implies $p\le\epsilon_m K$,
a contradiction.$\qquad\square$ \as

{\bf Proof of Lemma 4.1:}
The notation developed in the proof of Lemma 3.1 is used.

(i) Statement (7.2) implies $R_\Gamma \le R_\Gamma K/(n+m)$ since if
$Y_j\in C_{1,p,m}(X)$, then $Z_{0,j}$ has at most $R_\Gamma$ nonzero
elements. This implies $R_\Gamma = 0$, otherwise $n+1 \le n+m \le K$.
The remainder of the proof is similar to the proof of Lemma 3.2.ii.

(ii) Since $u$ is non-increasing, ${\rm trace}\, V_{n+m}(Z) \le B u(0)$.

(iii) If ${\rm trace}\, V_{n+m}(Z)$ is not bounded away from zero, then
there exists a sequence $V_{n+m}(Z_j) \to 0$. This implies
$p = \ave \{\psi(\bz' V_{n+m}^{-1}(Z_j)\bz)\} \to K > p$,
a contradiction.

(iv) If ${\rm trace}\, V_{n+m}^{-1}(Z)$ is not bounded above, then (7.3) holds.
This contradicts the condition on $\epsilon_m$ unless $R_\Gamma = 0$. However,
if $R_\Gamma = 0$ and $\alpha_j \to\infty$, then $p = 
\ave \{\psi(\alpha_j\bz'\Gamma_j\bz)\} \to K > p$, a contradiction.

(v)  The same argument used for (i) implies $R_\Gamma = 0$, then the same
argument used for (iii) implies ${\rm trace}\, V_{n+m}^{-1}(Z)$ must be
bounded.$\qquad\square$ \as

{\bf Proof of Theorem 5.1:}
Lemma 3.1.i insures the existence and uniqueness of
$V_{n+m}^{-1}(Z_k)$. Statement (5.2) implies a sequence $J$ exists such that
for $j\in J$, $A_n = V_{n+m}^{-1}(Z_{r(j)}) \to A > 0$. Arguments similar to
those used in the proof of Theorem 3.1 give 
$1\le\epsilon_m \psi(r(j)\btheta' A_j\btheta)$.
Unless $\btheta' A\btheta = 0$, $\psi(r(j)\btheta' A_j\btheta)\to A$,
which contradicts (5.1). Thus, $\btheta'V_{n+m}^{-1}(Z_r)\btheta\to 0$
as $r\to\infty$.

The notation developed in the proof of Lemma 3.2 is used in the 
remainder of this proof. Note that $A = \alpha\Gamma$ with $o < \alpha <
\infty$.

By (2.2), $\epsilon_m \psi(r(j)\alpha_j\btheta'\Gamma_j\btheta) \to c$
where $c = p - (n+m)^{-1} \sum_{1\le i\le n}\psi(\alpha\bx'_i\Gamma\bx_i)$.
The index set $J$ can be chosen so that 
$\Gamma_j^{1/2}\btheta/(\btheta'\Gamma_j\btheta)^{1/2} \to \bphi$ with
$\bphi'\bphi = 1$. Taking the limit in (7.1) over $j\in J$ and recalling
$P_\Gamma \Gamma^{1/2} = 0$ implies $P_\Gamma = cP_\Gamma \bphi\bphi'P_\Gamma$
and hence ${\rm rank}\, (P_\Gamma) = 1$ or ${\rm rank}\, (\Gamma) = p-1$.
Thus, for $\ba$ not proportional to $\btheta$, \mbox{$\ba' V_{n+m}^{-1}(Z_r)\ba$
is bounded away from zero. } \as

{\bf Proof of Theorem 5.2:}
The notation developed in the proof of Lemma 3.2 is used.
By (5.3), the subsequence $J$ can be chosen so that $V_{n+m}(Z_j) \to V$,
a nonzero positive semi-definite matrix. The matrix $V$ must be singular,
otherwise since $u$ is continuous Lemma 2.1.ii is contradicted when the
limit is taken. This implies $\alpha_j \to\infty$. Since
$\bx'_1 \Gamma_j\bx_1 
= \ave \{u(\alpha_j\bz'\Gamma_j\bz)\balpha_j(\bx'_1\Gamma_j\bz)^2\}$,
where the average is over $z\in Z_j$, if $\Gamma\bx_1 \ne 0$ then taking
the limit in the above statement gives $\bx'_1\Gamma\bx_1 \ge
(m+1)K(\bx'_1\Gamma\bx_1)^2/(n+m)$ or $n+m \ge (m+1)K$. This contradicts
the lower bound on $\epsilon_m$ and thus $\Gamma\bx_1 = 0$. The upper bound
on $\epsilon_m$ and (7.3) imply ${\rm rank}\, (\Gamma) = p-1$, and so since
$X$ is in general position about the center, $\Gamma\bx_i\ne 0$ for $i\ne 1$.
This implies $\alpha_j\bx'_i\Gamma_j\bx_i \to\infty$ or
$u(\alpha_j\bx'_i\Gamma_j\bx_i) \to 0$ for $1\ne 1$, and thus if
$\ba'\bx_1 = 0$, then $\ba' V_{n+m}(Z_j)\ba \to \ba' V\ba = 0$.
Since $V$ is  nonzero, $\bx'_i V\bx_1 > 0$.
The theorem follows since the arguments can be applied to any
convergence subsequence of $\ba' V_{n+m}(Z_k)\ba$ for $\ba'\bx_1 = 0$ 
or for \mbox{$\ba = \bx_1$.$\qquad\square$}

\medskip

\centerline{\bf R{\footnotesize{EFERENCES}}}
\medskip

Devlin, S.J., Gnanadesikan, R.\ and Kettenring, J.R. (1981).
Robust estimation of dispersion matrices and principle components.
{\it JASA} {\bf 76}, 354-362.

Donoho, D.L. (1982).  Breakdown properties of multivariate location
estimators. Ph.D.\ Qualifying Paper, Department of Statistics,
Harvard University.

Donoho, D.L. and Huber, P.J. (1983).
The notion of breakdown point. In {\it Festschrift in Honor of Erich
Lehmannn}, K.\ Doksum and J.L.\ Hodges, eds., Wadsworth, Belmont CA.

Hampel, F.R. (1971).
A general qualitative definition of robustness.
{\it Ann.\ Math.\ Statist.} {\bf 42}, 1887-1896.

Huber, P.J. (1981).
{\it Robust Statistics}. Wiley, New York.

Maronna, R.A. (1976).
Robust $M$-estimators of multivariate location and scatter.
{\it Ann.\ Stat.} {bf 4}, 51-67.

Stahel, W.A. (1981).
{\it Breakdown of covariance estimators.}
Research Rept.\ No.\ 31, E.T.H., Zurich.

Tyler, D.E. (1985).
Existence and uniqueness of the $M$-estimators of multivariate
location and scatter.

\end{document}